\documentclass[conference]{IEEEtran}
\IEEEoverridecommandlockouts

\usepackage[numbers,sort]{natbib} 
\usepackage{mathtools}
\mathtoolsset{showonlyrefs=true} 
\usepackage{amsmath,amssymb,amsfonts}
\usepackage{algorithmic}
\usepackage{graphicx}
\usepackage{textcomp}
\usepackage{xcolor}
\def\BibTeX{{\rm B\kern-.05em{\sc i\kern-.025em b}\kern-.08em
    T\kern-.1667em\lower.7ex\hbox{E}\kern-.125emX}}

\newcommand{\be}{\begin{equation}}
\newcommand{\ee}{\end{equation}}
\newcommand{\prob}[1]{\hbox{Pr}\left\{#1\right\}}

\newcommand{\kk}{\mathbf{k}}
\newcommand{\vv}{\mathbf{v}}
\newcommand{\pp}{\mathbf{p}}
\newcommand{\KK}{\mathbf{K}}

\newcommand{\eps}{\varepsilon}

\newcommand{\bb}[1]{\mathbf{#1}}

\begin{document}

\title{Bounding Means of Discrete Distributions\\
}

\author{\IEEEauthorblockN{Eric Bax}
\IEEEauthorblockA{\
\textit{Verizon}\\
Playa Vista, CA, USA \\
ebax@verizonmedia.com}
\and
\IEEEauthorblockN{Fr\'ed\'eric Ouimet}
\IEEEauthorblockA{\
\textit{McGill University}\\
Montreal, QC, Canada \\
frederic.ouimet2@mcgill.ca}}

\maketitle

\begin{abstract}
We introduce methods to bound the mean of a discrete distribution (or finite population) based on sample data, for random variables with a known set of possible values. In particular, the methods can be applied to categorical data with known category-based values. For small sample sizes, we show how to leverage knowledge of the set of possible values to compute bounds that are stronger than for general random variables such as standard concentration inequalities.
\end{abstract}


\section{Introduction}

In big-data settings, the data often comes from a sample of a population of interest, and we want to infer information about the population from our data. In this paper, we explore methods to bound population averages using a sample of observations with known category-based values. The observations are assumed to be independent and identically distributed (i.i.d.) from an unknown discrete distribution, or, equivalently, randomly drawn with replacement from our target population, with the goal of bounding the average value (or expectation) of the distribution.

For example, data on subscribers for a new service may include which tier of service each subscriber selected, and subscription pricing may be based on the tier of service. Based on this data, we may want to bound the average revenue per additional subscriber, under the assumption that new subscribers will select tiers of service in the same manner as present subscribers. Similarly, we may want averages conditioned on how the subscribers discovered the service, subscriber location and demographics, and other information, leading us to treat our big data set as many small data sets, each requiring estimated statistics.

As another example, suppose we have data for a sample of our subscribers, collected via reservoir sampling \cite{vitter84,vitter85,olken90,olken93,tille06} (with replacement). From that data, we may want to bound the averages of some statistics over all subscribers, such as the number of visits to a website or the number of uses of a service per week. In this case, the numbers categorize the subscribers as well as being the values of interest. (The motivation for sampling may be efficiency, preservation of privacy \cite{denning80,domadiya19}, or both.)

One method to obtain bounds on population or distribution averages is to apply concentration inequalities to the category values of the observations. For example, we can use bounds on the sample average that only rely on the range of the category values, such as \citet{hoeffding63}, \citet{azuma67}, and \citet{mcdiarmid89} bounds, or use bounds that also rely on the sample variance, sometimes called empirical Bernstein bounds \cite{bernstein37,audibert08,maurer09}. Since we know the category value of each observation in our sample though, we have some extra information that those bounds do not utilize. It is therefore reasonable to hope that we can exploit that information to produce stronger bounds.

Instead of concentration inequalities, our general technique is validation by inference \cite{dunn61,bax_voting,bax_val_by_inference}: first find a set of distributions that includes all those for which our observations are likely (in the sense of not being too far out in the tails of those distributions), then identify distributions in that ``likely set'' that have minimum or maximum means. The minimum and maximum are lower and upper bounds, respectively, on the mean of the distribution that generated our sample, with probability of bound failure no more than the probability that our sampled observations are ``too far'' out in the tail of their distribution. The bounds in this paper, both from concentration inequalities and from validation by inference, are PAC (probably approximately correct) bounds \cite{valiant84}.

The next section reviews some useful mathematical tools. Section~\ref{section_bb} explains how to perform validation by inference with a likely set based on treating each category's sample count as a binomial random variable and using simultaneous bounds over categories to infer bounds on distribution averages. Section~\ref{section_nb}, which is the main contribution of this paper, shows how to improve the validation by inference method by treating combined category sample counts as binomial samples, with a nested pattern of combinations. Section~\ref{section_test} compares those methods to using known concentration inequalities. Section~\ref{section_future} outlines potential future work, including forming a likely set based on a multivariate normal approximation to the multinomial distribution that generated the vector of category sample counts.

\section{Preliminaries}
Assume that each observation can be placed into one of $m$ possible categories, and let $\KK = (K_1, \ldots, K_m)$ denote the vector of sample counts for the $m$ categories. (The observed sample counts are not capitalized, i.e., $\kk = (k_1,\dots,k_m)$). Let $n = K_1 + \ldots + K_m$ be the total number of observations ($n$ is fixed and known). Given that the observations are i.i.d., we have $\KK\sim \mathrm{Multinomial}\hspace{0.2mm}(n,\pp^{\star} = (p^{\star}_1, \ldots, p^{\star}_m))$, where $p_i^{\star}$ denotes the probability of any observation to fall in the $i$-th category. Let $\vv = (v_1, \ldots, v_m)$ denote the vector of values assigned to the $m$ categories, and assume that $v_1 < \ldots < v_m$. (Combine any categories that have equal values.) Our goal is to compute PAC bounds \cite{valiant84} on $\pp^{\star} \cdot \vv$ (``$\cdot$'' is the dot product), the out-of-sample expectation of category value, with some specified probability of bound failure at most $\delta > 0$.

If there are two categories ($m = 2$), then the sample count for the second category is a binomial random variable, so we can compute a bound using binomial inversion \cite{hoel54,langford05}, as follows. For any $p\in [0,1]$, define $B(n, k, p)$ to be the left tail (the cdf) of the binomial distribution:
\be
B(n,k,p) = \sum_{i=0}^{k} \binom{n}{i} p^i (1-p)^{n-i}, \quad k\in \{0,1,\dots,n\}. \label{left_tail}
\ee
Then, with probability at least $1 - \delta$, the binomial inversion upper bound:
\be
p_{+}(n, k, \delta) = \max \{p\in [0,1] : B(n, k, p) \geq \delta\}. \label{bin_inv_upper}
\ee
is at least the probability of an event that occurs $k$ times in $n$ independent Bernoulli trials. By definition, this bound is sharp, in the sense that the bound failure probability is $\delta$. Also, it is easy to compute, because $B(n, k, p) \geq \delta$ for all $p \leq p_{+}(n, k, p)$ and $B(n, k, p) < \delta$ for all $p > p_{+}(n, k, p)$. So we can use binary search over $p \in [0, 1]$, with precision $\frac{1}{2^s}$ after $s$ search steps.

So, for $m = 2$, we have, with probability at least $1 - \delta$,
\be
p^{\star}_2 \leq p_{+}(n, k_2, \delta).
\ee
Thus, we deduce
\be
\pp^{\star} \cdot \vv \leq [1 -  p_{+}(n, k_2, \delta)] v_1 + p_{+}(n, k_2, \delta) v_2,
\ee
since increasing $p^{\star}_2$ increases $\pp^{\star} \cdot \vv$. (Remember $v_1 < v_2$.)

Similarly, for a lower bound on $\pp^{\star} \cdot \vv$, we have, with probability at least $1 - \delta$,
\be
p^{\star}_1 \leq p_{+}(n, k_1, \delta),
\ee
and subsequently,
\be
\pp^{\star} \cdot \vv \geq p_{+}(n, k_1, \delta) v_1 + [1 - p_{+}(n, k_1, \delta)] v_2.
\ee

Note that computing the lower bound is equivalent to reversing the order of elements in $\kk$ and $\vv$, so that the values are in descending order, then applying the procedure used to compute the upper bound. For two-sided (simultaneous upper and lower) bounds, use $\frac{\delta}{2}$ in place of $\delta$ in each bound. This is called the Bonferroni correction \cite{bonferroni36} or the union bound, since the probability of the union of events is at most the sum of the probabilities of events. (In this case, the two events are bound failures for the upper bound and for the lower bound.)

For $m > 2$, $\pp^{\star}$ specifies a multinomial distribution, and computing a bound is more challenging. For any probability vector $\pp$, define the random vector $\KK_{\pp}\sim \mathrm{Multinomial}\hspace{0.2mm}(n,\pp)$ and, for a vector of sample counts $\bb{k}\in \mathbb{N}_0^m \cap n \mathcal{S}_m$, define
\be
L^{\star}(\bb{k}) = \{\pp\in \mathcal{S}_m : \prob{\bb{K}_{\bb{p}} \cdot \bb{v} \leq \bb{k} \cdot \bb{v}} \geq \delta\},
\ee
where $\mathcal{S}_m = \{\bb{s}\in [0,1]^m : \sum_{i=1}^m s_i = 1\}$.
This is called the \textit{likely set}, namely the set of probability vectors $\pp$ that are likely to have generated our sample in the sense that the observed average $\kk \cdot \vv$ is not in the left $\delta$-tail of the distribution of $\KK_{\pp} \cdot \vv$.
Then an upper bound that holds with probability at least $1 - \delta$ is
\be
\pp^{\star} \cdot \vv \leq \max_{\pp \in L^{\star}} \pp \cdot \vv.
\ee

We do not know how to compute this ideal upper bound directly. (Unlike for $m = 2$, we cannot use binary search for $m > 2$.) Therefore, we must use bounds based on alternative likely sets rather than $L^{\star}$. For a valid bound, we require that the (unknown) probability vector $\pp^{\star}$ that generated the observations be in the likely set with probability at least $1 - \delta$. 


As an aside, it is computationally feasible to determine whether a given $\pp$ is in $L^*$, by using Monte Carlo sampling \cite{hope68} or an ``exact test'' \cite{keich06,jann08}. (The referenced methods test based on likelihood or Pearson's $\chi^2$ \cite{pearson00,siotani84,cressie84,read88,balakrishnan18}, but they can be easily adapted to our $L^*$.). Nonetheless, identifying a $\pp$ in $L^*$ that maximizes $\pp \cdot \vv$ is more challenging, because it requires optimization. 



The likely sets we use in this paper are convex sets with linear boundaries, making it possible to maximize $\pp \cdot \vv$ using linear programming. However, because of the specific sets we use, we can apply simpler optimization methods than needed for general linear programs.

\section{Bonferroni Box} \label{section_bb}
Each category's sample count has a binomial distribution. Therefore, we can use binomial inversion to bound the probability of membership in each category. In this section, we combine binomial inversions for each category using a Bonferroni correction, forming a likely set that is a rectangular prism (a ``box'') that contains the generating distribution $\pp$ with probability at least $1 - \delta$. Given such a simple shape for the likely set, it is easy to find the probability vector $\pp$ in the set that maximizes $\pp \cdot \vv$ to produce an upper bound on the expectation of the category value.

Define a binomial inversion lower bound:
\be
p_{-}(n, k, \delta) = \min \{p\in [0,1] : 1 - B(n, k - 1, p) \geq \delta\}.
\ee
and use it to define a \textit{Bonferroni box}:
\be
\textstyle L_B = \text{\Large $\times_{i\in [m]}$}  \left[p_{-}(n, k_i, \tfrac{\delta}{2 m}), p_{+}(n, k_i, \tfrac{\delta}{2 m})\right],
\ee
where $[m] = \{1, \ldots, m\}$. Then, with probability at least $1 - \delta$, $\pp^{\star} \in L_B$, since
\be
\prob{p^{\star}_i \not\in [p_{-}(n, k_i, \tfrac{\delta}{2 m}), p_{+}(n, k_i, \tfrac{\delta}{2 m})]} \leq \tfrac{\delta}{m}, ~\forall i \in [m],
\ee
which implies
\be
\prob{\exists i \in [m] : p^{\star}_i \not\in [p_{-}(n, k_i, \tfrac{\delta}{2 m}), p_{+}(n, k_i, \tfrac{\delta}{2 m})]} \leq \delta,
\ee
based on the Bonferroni correction / union bound. Hence, the maximum
\be
\max_{\pp \in L_B} \pp \cdot \vv
\ee
is an upper bound on $\pp^{\star} \cdot \vv$, with probability at least $1 - \delta$.

To find the maximizing $\pp \in L_B$, first assign each $p_i$ to its lower bound. Call the difference between one and the sum of the $p_i$ values the headroom, and update it at each step. For each $p_i$ value starting with $p_m$ and working back to $p_1$, if the headroom is greater than zero, add the headroom or the difference between the upper and lower bound for $p_i$, whichever is least. This allocates the probability mass to the rightmost $p_i$ values, to the extent allowed by the upper bounds for the rightmost values while also allocating at least the lower bounds for the leftmost elements. Since $v_1 < \ldots < v_m$, this maximizes $\pp \cdot \vv$.

For the lower bound, just reverse $\kk$ and $\vv$ and apply the procedure we outlined for the upper bound. For simultaneous upper and lower bounds, there is no need to divide $\delta$ by 2, since both bounds use the same $2m$ single-category upper and lower bounds. (Since we infer both bounds from the same set of single-category ``basis'' bounds, we only need one union bound over the basis bounds.)


\section{Bonferroni Nest} \label{section_nb}
Rather than use a separate binomial inversion bound for the probability of each category, it is possible to apply binomial inversion bounds to the sums of probabilities of combinations of categories. Each bound on the probability of a combination of categories implies a constraint on a sum over the entries in $\pp$ that correspond to the categories. Together, these constraints imply a bound on $\pp^{\star} \cdot \vv$. The resulting bound can be an improvement over a bound based on computing separate binomial inversion bounds for each category. 

To see why, recall that the variance of a binomial distribution is $np(1-p)$, so the standard deviation of each category's sample count is $\sqrt{np(1-p)}$, where $p$ is the category probability. The differences between frequencies and binomial inversion bounds scale approximately with the standard deviation of the category's sample count, divided by the total number of observations, i.e.,
\be
\frac{\sqrt{np(1-p)}}{n} = \sqrt{\frac{p(1-p)}{n}}.
\ee
If we combine $c$ categories, each with probability $p$, then the combined probability is $cp$. So the difference between the resulting binomial inversion bound and the combined frequency scales as
\be
\sqrt{\frac{cp(1-cp)}{n}} \leq \sqrt{c} \sqrt{\frac{p(1-p)}{n}}.
\ee
This is about $\sqrt{c}$ times the difference between the frequency and the bound for a single category. In contrast, if we bound $c$ categories separately and sum the bounds, then the difference between the sum of frequencies and the sum of bounds is $c$ times the difference for a single category. Therefore, we can get tighter bounds on combined categories by summing frequencies then bounding instead of bounding frequencies then summing.

Let
\be
t_0 = 0,
\ee
\be
\textstyle \forall i \in [m-1] : t_i = p_{-}(n, \sum_{j=1}^{i} k_j, \frac{\delta}{m - 1}),
\ee
and
\be
t_m = 1.
\ee
Then each $t_i$ is a lower bound on $p_1^{\star} + \ldots + p_i^{\star}$, and the bounds hold simultaneously with probability at least $1 - \delta$. (The bound $t_m = 1$ follows from $\pp^{\star}$ being a probability vector.)

Let $L_N$ be the set of probability vectors $\pp$ that satisfy the lower-bound constraints:
\be
\forall i \in [m] : p_1 + \ldots + p_i \geq t_i.
\ee
Recall that $v_1 < \ldots < v_m$. So to maximize $\pp \cdot \vv$ over $\pp \in L_N$, we will place as little probability in earlier $p_i$ values, and as much in later ones, as possible.

Begin with $p_1$. Since $p_1 \geq t_1$, $t_1$ is the least probability that we can assign to $p_1$. So set $p_1 = t_1 = t_1 - t_0$.

For $i > 1$, we have the lower bound $p_1 + \ldots + p_i \geq t_i$. But the previous lower bound, $p_1 + \ldots + p_{i-1} \geq t_{i-1}$, forces us to assign at least $t_{i-1}$ in total to $p_1 + \ldots + p_{i-1}$. That leaves $t_i - t_{i-1}$ as the most of $t_i$ that we can assign to $p_i$ while assigning the minimum possible ($t_i$) to the sum $p_1 + \ldots + p_i $. (Assigning that minimum leaves as much probability as possible for $p_{i+1} + \ldots + p_{m}$.)

So assign each $p_i = t_i - t_{i-1}$ to maximize $\pp \cdot \vv$. The resulting $\pp$ is a probability vector: the sum is one since $t_0 = 0$ and $t_m = 1$, and each entry is nonnegative since $t_0 \leq \ldots \leq t_m$. (Binomial inversion bounds increase monotonically in $k$).

For a lower bound, follow the same procedure, but reverse $\kk$ and $\vv$. For simultaneous upper and lower bounds, use $\frac{\delta}{2}$ in place of $\delta$ for each bound, because the nested bounds nest in different directions (right and left in the original category ordering) in the two bounds, collecting different sets of categories. This is not a disadvantage compared to the Bonferroni box bound, because that uses $2m$ individual bounds for simultaneous upper and lower bounds, and this Bonferroni nest bound uses only $2(m - 1)$.

\section{Comparisons} \label{section_test}

\begin{figure}[htbp]
\centerline{\includegraphics{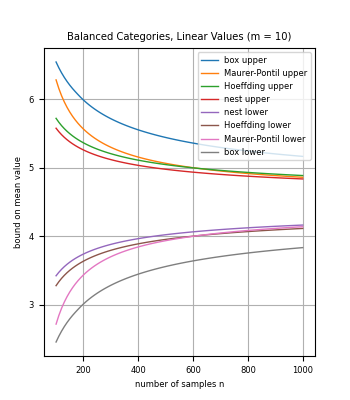}}
\caption{All categories have the same sample counts, and category values increase linearly from 0 to 9. For $n = 100$ to 1000 observations, Bonferroni nest bounds are the tightest, and Bonferroni box bounds are the least tight. Maurer-Pontil bounds outperform Hoeffding bounds for more than about 600 observations; they also eventually outperform Bonferroni nest bounds, starting at a few thousand observations (not shown).}
\label{bal_lin_m10}
\end{figure}

\begin{figure}[htbp]
\centerline{\includegraphics{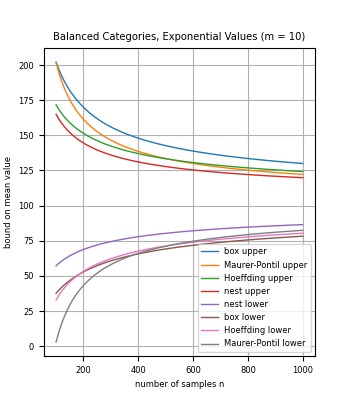}}
\caption{All categories have the same number of observations. Category values increase exponentially from $2^0$ to $2^9$, making later categories play a stronger role in determining the average. Bonferroni nest bounds are the tightest for up to 1000 observations. Maurer-Pontil bounds are the least tight for small samples, but are the tightest for several thousand or more observations (not shown).}
\label{bal_exp_m10}
\end{figure}

\begin{figure}[htbp]
\centerline{\includegraphics{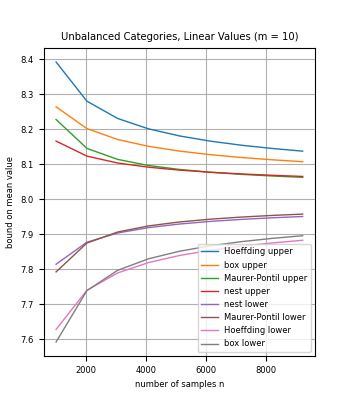}}
\caption{Sample count per category increases exponentially. Bonferroni nest bounds are the tightest for up to a few thousand observations. Then Maurer-Pontil bounds are best.}
\label{unb_lin_m10}
\end{figure}

\begin{figure}[htbp]
\centerline{\includegraphics{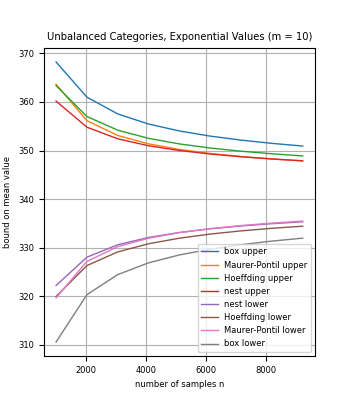}}
\caption{Number of observations doubles from category to category, and category values start with 1 and also double from category to category. As a result, the last category plays an outsize role in determining the average category value per observation. Bonferroni bounds are superior for fewer than several thousand observations.}
\label{unb_exp_m10}
\end{figure}

\begin{figure}[htbp]
\centerline{\includegraphics{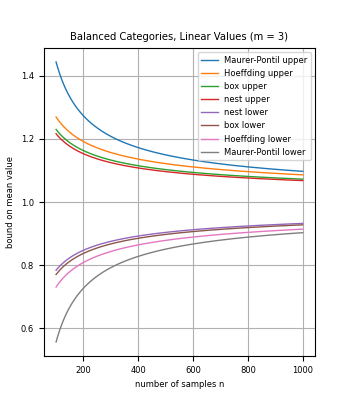}}
\caption{With $m=3$ categories, both Bonferroni methods outperform the single-variable methods. Bonferroni nest bounds are stronger than Bonferroni box bounds.}
\label{bal_lin_m3}
\end{figure}

\begin{figure}[htbp]
\centerline{\includegraphics{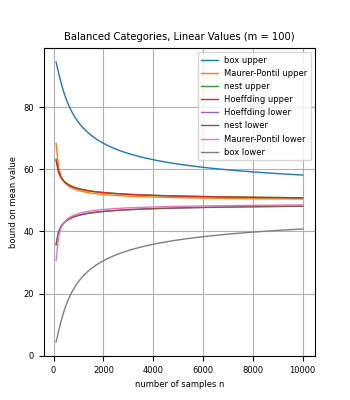}}
\caption{With $m = 100$ categories, the Bonferroni methods divide $\delta$ so much that they become less effective than the single-variable methods. Hoeffding bounds are the best for the smallest numbers of observations, then Maurer-Pontil bounds are best. For the smallest numbers of observations, the Bonferroni nest bound is between the Hoeffding and Maurer-Pontil bounds, and it remains close to them for large numbers of observations.}
\label{bal_lin_m100}
\end{figure}

We can think of our in-sample data as a set of $n = k_1 + \ldots + k_m$ observations, with $k_i$ observations of value $v_i$ for each $i$. The range is $r = v_m - v_1$, the sample mean is $\hat{\mu}_n = \kk \cdot \vv / n$, and the sample variance is $\hat{\sigma}_n^2 = \sum_i k_i (v_i - \hat{\mu}_n)^2 / (n - 1)$.

Then an empirical Bernstein bound, meaning a bound based on the sample variance, by Maurer and Pontil \cite{maurer09} has:
\be
\eps_{\mathrm{MP}} = \sqrt{\frac{2 \hat{\sigma}_n^2 \ln(\frac{4}{\delta})}{n}} + \frac{7 r \ln(\frac{4}{\delta})}{3(n - 1)}
\ee
for simultaneous lower and upper bounds $\hat{\mu}_n - \eps_{\mathrm{MP}}$ and $\hat{\mu}_n + \eps_{\mathrm{MP}}$.

The Hoeffding bound \cite{hoeffding63} has
\be
\eps_{\mathrm{H}} = r \sqrt{\frac{\ln(\frac{2}{\delta})}{2 n}}
\ee
for simultaneous lower and upper bounds $\hat{\mu}_n - \eps_{\mathrm{H}}$ and $\hat{\mu}_n + \eps_{\mathrm{H}}$.

The Hoeffding bound relies on worst-case assumptions about variance, making it simpler to implement and analyze than empirical Bernstein bounds. (Hoeffding's paper \cite{hoeffding63} does include bounds that incorporate variation, but the simpler bound presented here is more often used.) Conceptually, empirical Bernstein bounds must validate variance simultaneously with the mean, so they suffer for very small sample sizes, then tend to outperform worst-case variance bounds for distributions with limited variance.

In this section, we compare the Bonferroni box bound and the Bonferroni nest bound to the Hoeffding and Maurer and Pontil bounds using empirical tests. In general, the Bonferroni nest bound is superior to the Bonferroni box bound. The Bonferroni nest bound offers non-trivial bounds even for very small sample sizes ($n < 100$), and it tends to outperform the other bounds for small to moderate sample sizes. But for large sample sizes the Maurer and Pontil bound is best. In computer science, we are often concerned about what happens as the problem size grows. In applied statistics, we are often concerned with what happens as the sample size shrinks. When we mix the two, the right answer depends on the specifics of the problem.

Our initial tests involve two types of category sample count vectors $\kk$: \textit{balanced} have all entries equal, meaning that each category has an equal number of observations, and \textit{unbalanced} have each entry double the previous entry, meaning that higher-value categories have more observations. The tests also involve two types of value vectors $\vv$: \textit{linear} have values $0, 1, \ldots, m - 1$, and \textit{exponential} have values $2^0, 2^1, \ldots, 2^{m-1}$. To begin, we set $m = 10$ categories; we will vary the number of categories later. For all tests in this section, $\delta = 0.05$ is the probability of bound failure allowed for simultaneous two-sided bounds.

Figure~\ref{bal_lin_m10} compares Bonferroni box, Bonferroni nest, Hoeffding, and Maurer-Pontil bounds. Bonferroni nest bounds are the best for fewer than 1000 observations, though Maurer-Pontil bounds are best for a few thousand observations or more (not shown). Maurer-Pontil bounds are weaker than Hoeffding bounds for fewer than 600 observations -- since Maurer-Pontil bounds validate variance in order to validate the mean, they are in some sense performing more validations than Hoeffding bounds, so they require more observations to begin being effective. Bonferroni box bounds are the least effective, indicating that the gains from nesting validations make the Bonferroni approach better than Hoeffding and Maurer-Pontil bounds for small numbers of observations.

The tests for Figure~\ref{bal_exp_m10} are similar to those for Figure~\ref{bal_lin_m10}, except that category values increase exponentially instead of linearly, starting with 1 and doubling from category to category. Bonferroni nest bounds are the best for fewer than 1000 observations. As in the previous tests, Maurer-Pontil bounds are weaker than Hoeffding bounds for a few hundred observations, but surpass them starting at around 500 observations. Beyond a few thousand observations, Maurer-Pontil bounds are best.

The tests for Figure~\ref{unb_lin_m10} are similar to those for Figure~\ref{bal_lin_m10}, except that the number of observations increases exponentially from category to category, doubling from each category to the next. Bonferroni nest is the best upper bound for less than 4000 observations, and the best lower bound for less than 200 observations. Beyond those numbers, Maurer-Pontil bounds become the best.

The tests for Figure~\ref{unb_exp_m10} have both numbers of observations and category values increasing exponentially from category to category. Again, Bonferroni nest bounds are the best for fewer than several thousand observations. Then Maurer-Pontil bounds are stronger.

Now consider how the number of categories $m$ affects the bounds. For $m = 2$, the Bonferroni nest bound is a (sharp) binomial inversion bound. But as the number of categories increases, the Bonferroni correction divides $\delta$ into smaller $\delta$ values for the individual bounds. This is a disadvantage for the Bonferroni bounds compared to the single-distribution bounds: Hoeffding and Maurer-Pontil. Figure~\ref{bal_lin_m3} shows that for $m = 3$ categories, the Bonferroni methods are the most effective. (Compare to Figure~\ref{bal_lin_m10}, for $m = 10$.) Figure~\ref{bal_lin_m100} shows that for $m = 100$ categories, the best single-variable bound is better than the best Bonferroni bound, though the Bonferroni nest bound is quite close to the single-variable bounds.

\section{Merging Categories}
To avoid over-partitioning $\delta$ for the Bonferroni bounds, we can combine some neighboring categories, summing their sample counts and using their maximum category value as the value for the combined category, for upper bounds. (Use the minimum for lower bounds.) Combining categories decreases the number of categories $m$, so it partitions $\delta$ less, which helps strengthen bounds. But it also attributes all observations for the combined category to the worst-case value among those for the categories combined. That weakens bounds.

To minimize the worst-case weakening due to collecting observations from multiple categories and giving them the worst-case value, we want to minimize the maximum of the range of category values for any set of categories combined into a single category. We can do this via dynamic programming. Refer to a sequence of categories to be merged into a single category as a cluster. Refer to the difference between the maximum and minimum values for categories in a cluster as the cluster range. Let $c_{hj}$ be the minimum, over ways to put categories 1 to $j$ into $h$ clusters, of the maximum cluster range. Since merging categories 1 to $j$ into a single cluster gives cluster range $v_j - v_1$:
\be
\forall j \in [m] : c_{1j} = v_j - v_1.
\ee
For multiple clusters:
\be
\forall h \geq 2, j \geq h: c_{hj} = \min_{h \leq i < j} \max(c_{h - 1, i}, v_{j} - v_{i + 1})
\ee
because the last cluster can be categories $i + 1$ to $j$ for any $h \leq i < j$, making the last cluster range $v_j - v_{i + 1}$ and the (best-case) maximum of the other cluster ranges $c_{h - 1, i}$. Implementing this recurrence as a dynamic program with choice recovery produces a clustering with $h$ clusters, with minimum (over clusterings) maximum cluster range $c_{hm}$.

\begin{figure}[htbp]
\centerline{\includegraphics{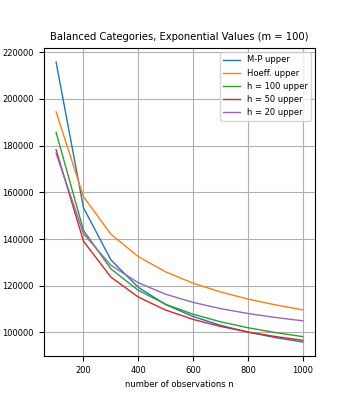}}
\caption{With $m = 100$ categories and exponentially increasing category values ($v_i = 2^{20 (i - 1) / m}$), merging categories can improve nested Bonferroni bounds. The number of categories after merging is $h$, and $h = 100$ is the nested Bonferroni bound with the original 100 categories. To show detail, only the upper bounds are plotted. For 100 observations, 20 merged categories outperforms the other bounds, but using 50 merged categories is superior for more observations. For 1000 or more observations, Maurer-Pontil bounds are best.}
\label{cluster_fig}
\end{figure}

For each cluster, sum the sample counts $k_i$ and take the worst-case of the values $v_i$ (maximum for upper bound, minimum for lower bound), to form a clustered sample count vector $\kk$ and value vector $\vv$, then apply a Bonferroni bound, and it will also apply to the original sample and values. Figure~\ref{cluster_fig} compares Maurer-Pontil and Hoeffding bounds to bounds based on merging categories before applying a nested Bonferroni bound. To show detail, only the upper bounds from two-sided bounds are shown, and the number of observations ranges from 100 to 800 (with equal numbers per category).

The category values are selected to increase (strongly) exponentially ($v_i = 2^{20 (i - 1) / m}$) in order to illustrate the possibility for merging categories to improve bounds. (Reducing the value 20 in the exponent reduces the range of observation counts for which any other bound shown is superior to the best of the Maurer-Pontil and Hoeffding bounds. Setting the value to one leaves no such range.)

For $h = 100$ merged categories, every original category is a merged category, so the bound is the original nested Bonferroni bound. For $h = 50$, the dynamic program merges categories 1 through forty into a single category, merges the next five categories together, the next three, the next two pairs, and then keeps each remaining category as a separate category in the clustering. For $h = 20$, the first 70 categories are merged, then five, three, and three pairs are merged, leaving the last 15 as separate categories.

\section{Nearly Uniform Category Bounds}
The Bonferroni nest upper bound on $\pp^{\star} \cdot \vv$ is based on $m - 1$ simultaneous frequency bounds for sets of categories, where $m$ is the number of categories. Each frequency bound uses $\frac{\delta}{m - 1}$ for $\delta$, so that the probability of any frequency bound failing is at most $\delta$. As the number of categories grows, $\delta$ becomes so partitioned that the binomial inversion bounds become ineffective. Using nearly uniform bounds \cite{bax_compression,bax_avg} counters this problem by allowing some frequency bound failures in exchange for an increase in $\delta$ for each frequency bound. To understand how this works, consider the probability that at least five events out of 100 different events occur, assuming that each event has at most a one percent probability of occurrence. The worst-case joint distribution is that any one event occurring implies that exactly four others occur as well, so either no events occur or five occur together. So each one percent probability of five events occurring together requires five percent from the sum of the probabilities of the events. With 100 events each with one percent probability, then, the maximum probability of five events occurring together is at most $(100)(0.01)/5$, which is 20\%.

In our case, if each of $m - 1$ events has probability at most $(a + 1)\frac{\delta}{m - 1}$, then the probability of $a + 1$ or more events is at most
\be
\frac{(m - 1) (a + 1) \frac{\delta}{m - 1}}{a + 1} = \delta.
\ee
So allowing $a$ bound failures allows us to use $(a + 1)\frac{\delta}{m - 1}$ for $\delta$ in the frequency bounds and still have bound failure probability at most $\delta$ for the Bonferroni nest upper bound. But we must adjust the bound to account for up to $a$ frequency bound failures.

The Bonferroni nest upper bound process produces a probability vector $\pp$ that maximizes $\pp \cdot \vv$ under the frequency bound constraints. Each frequency bound failure allows the probability mass from some $p_{i - 1}$ to move to $p_i$, increasing the bound by $p_{i - 1} (v_i - v_{i - 1})$. Sequential bound failures allow probability mass from sequential entries to all accumulate in the entry to the right of the sequence. So $h$ bound failures in a row, just prior to the bound on $p_1 + \ldots + p_i$, increases the nest bound by
\be
\Delta(h, i) = \sum_{b = 1}^{h} p_{i-b} (v_i - v_{i - b}).
\ee

We can use dynamic programming to maximize these increases for $a$ allowed bound failures. Let $c_{ij}$ be the maximum nest bound increase for $j$ bound failures before the bound for $p_1 + \ldots + p_i$ and no bound failures after. The base cases are for all bounds failing before the bound for $p_1 + \ldots + p_i$, so $i = j + 1$:
\be
\forall i \in \{1, \ldots, a\}: c_{i, i - 1} = \Delta(i - 1, i).
\ee
The general recurrence maximizes over all possible numbers of bound failures $h$ immediately preceding the bound for $p_1 + \ldots + p_i$, combining the effect of those bound failures with the effect of previous ones:
\be
\forall j \in \{1, \ldots, a\}, i \in \{j + 2, \ldots, m\}:
\ee
\be
c_{ij} = \max_{0 \leq h \leq j} \left[c_{i - 1 - h, j - h} + \Delta(h, i) \right].
\ee
Add $c_{m, a}$ to the Bonferroni nest bound for a nearly uniform Bonferroni nest bound.

For a nearly uniform Bonferroni nest lower bound, reverse $\kk$ and apply the Bonferroni nest upper bound process to get $\pp$, as we did for the Bonferroni nest lower bound. Then apply the dynamic program to those $\pp$ values with $\vv$ reversed and multiplied by negative one, and subtract the resulting $c_{m, a}$ from the Bonferroni nest lower bound. For simultaneous lower and upper bounds, use $(a + 1) \frac{\delta}{2}$ in place of $\delta$ in each of the lower and upper bounds, and apply the corrections for allowing $a$ bound failures.

\begin{figure}[htbp]
\centerline{\includegraphics{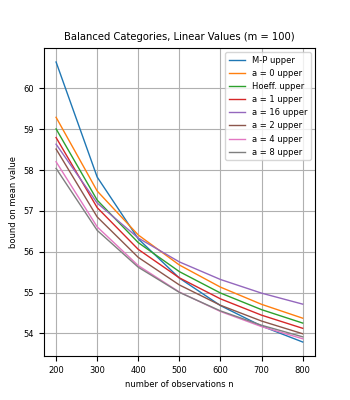}}
\caption{With $m = 100$ categories and linear category values ($v_i = i - 1$), nearly uniform Bonferroni nest bounds outperform the original Bonferroni nest bound for small observation counts. Only upper bounds are plotted, to show detail. The number of bound errors allowed in the nearly uniform bounds is $a$. So $a = 0$ is the original Bonferroni nest bound. Allowing eight errors outperforms allowing fewer for a few hundred observations, but there is crossover before 800 observations. Maurer-Pontil bounds are the best among those shown for 800 or more observations.}
\label{nu_fig}
\end{figure}

Figure~\ref{nu_fig} shows that using a nearly uniform technique can improve Bonferroni nest bounds to make them better than either Hoeffding or Maurer-Pontil bounds for small numbers of observations. To show detail, only upper bounds are shown for two-sided bounds. As in Figure~\ref{bal_lin_m100}, there are equal numbers of observations per category for 100 categories, with linearly increasing category values. The number of allowed bound failures is $a$, so $a = 0$ is the original Bonferroni nest bound.

For each number of observations, the better of Maurer-Pontil and Hoeffding bounds is superior to the original Bonferroni nest bound. But allowing a few errors produces stronger bounds for small numbers of observations. Allowing too many errors weakens bounds relative to the original Bonferroni nest bound as the number of observations increases, as indicated by the results for $a = 16$.

\section{Discussion} \label{section_future}
We have shown that we can use the information that a distribution has support over a discrete set of values to produce mean bounds that are stronger than those from standard methods for general distributions, for small numbers of observations. The advantage offered by our method, Bonferroni nest bounds, generally decreases as the number of values in the support of the distribution increases. Informally, more values mean that the distribution is ``less discrete'', in the sense that it is a closer approximation to a continuous distribution. For Bonferroni nest bounds, the problem with more values is that they require more simultaneous probability bounds to produce a mean bound.

We explored two methods to address this issue: combining sample counts for neighboring values into counts for the worst of those values, and allowing some probability bound failures. Both methods improve bounds for some distributions, and both have their limitations. In the future, it would also be interesting to consider partitioning the bound failure probability, $\delta$, in some other way than equally over the probability bounds.

It would also be interesting to develop a general mean bound for discrete-valued distributions, using the set of values as an input to produce a worst-case optimal Bonferroni nest bound with a combination of merged values, nearly uniform bounds, and partitioning of $\delta$. It is trivial to combine these tactics, but it may not be trivial to optimize over combinations.

For moderate and large sample sizes, it would be useful to incorporate information about the sample counts for the values, into bound selection, generalizing the technique of empirical Bernstein bounds. It may be possible to improve Maurer-Pontil bounds for discrete-valued distributions by deriving a tighter bound on the variance of the distribution than the one for general distributions.

It may also be possible to derive mean bounds for discrete-valued distributions based on concentration inequalities for multinomial distributions. It is well-known that the multinomial distribution has an asymptotic multivariate normal approximation \cite{pearson00,cressie84,siotani84,read88,gaunt16,ouimet21}. A likely set based on a multivariate normal could be an ellipsoid, making it easy to maximize $\pp \cdot \vv$ using quadratic programming. Developing such a likely set has a few challenges: we would need effective non-asymptotic ellipsoid probability bounds (rather than just asymptotic approximations), and we would need to find a reasonably-sized ellipsoid that contains all likely generating distributions $\pp$, instead of just showing that a reasonably-sized ellipsoid has most of the probability mass for a single distribution. There has been some progress toward these goals \cite{ouimet21}, but not a complete solution.

Finally, it may be interesting to apply methods for discrete-valued distributions to general distributions, which may have support over continuous sets. It is possible to merge all observed values into a discrete set of subranges, which may each be continuous, and use the value from the subrange that gives the worst-case mean -- largest value for an upper bound, smallest for a lower bound -- with the subrange frequencies as category values and frequencies in the methods from this paper. For some distributions, this may yield effective mean bounds.

\section*{Acknowledgments}

F.\ Ouimet is supported by postdoctoral fellowships from the NSERC (PDF) and the FRQNT (B3X supplement and B3XR).

\bibliographystyle{IEEEtranN}
\bibliography{bax}

\end{document}